\documentclass[a4paper, 11pt]{article}
\usepackage{mathrsfs}
\usepackage{latexsym,epsfig,amsmath,bm}
\usepackage{amssymb}
\usepackage{ebezier}
\usepackage{epic,eepic}
\usepackage{float}
\newcounter{prostredi}
\def\theprostredi{\arabic{prostredi}}

\makeatletter
\def\@begintheorem#1#2{\trivlist
   \item[\hskip \labelsep{\bfseries #1\ #2.}]} 
\def\@opargbegintheorem#1#2#3{\trivlist
      \item[\hskip \labelsep{\bfseries #1\ #2.\ (#3.)}]} 
\def\@endtheorem{\endtrivlist}
\makeatother

\makeatletter \@addtoreset{equation}{section}

\makeatother

\newtheorem{claim}{Claim}

\newcommand{\qed}{{$\qquad\square$\bigbreak}}

\newtheorem{theorem}{Theorem}

\newtheorem{corollary}[theorem]{Corollary}

\title{\huge Forbidden pairs for equality of edge-connectivity and minimum degree}

\author{\emph{Junfeng Du}\\ [0.3cm]
\small School of Mathematics and Statistics\\
\small Beijing Institute of Technology\\
\small Beijing 100081, P.R.\ of China\\
\small e-mail: djfdjf1990@163.com\\
\and \emph{Ziwen Huang}\\ [0.3cm]
\small School of Mathematics and Computer Science \\
\small \& Center of Applied Mathematics\\
\small   Yichun University\\
\small  Yichun 336000, P.R.\ of China\\
\small  e-mail: zwhuang@aliyun.com\\
\and \emph{Liming Xiong}\\ [0.3cm]
\small  School of Mathematics and Statistics\\
\small  Beijing Institute of Technology\\
\small  Beijing 100081, P.R.\ of China\\
\small  e-mail: lmxiong@bit.edu.cn\\}

\begin{document}
\date{}
\maketitle

\begin{abstract}
Let $\mathcal{H}$ be a class of given graphs. A graph $G$ is said to
be $\mathcal{H}$-free if $G$ contains no induced copies of $H$ for
any $H \in \mathcal{H}$. In this article, we characterize all pairs
$\{R,S\}$ of graphs such that every connected $\{R,S\}$-free graph
has the same edge-connectivity and minimum degree.

 \vspace*{0.2cm}\noindent
  {\bf Keywords}: forbidden subgraph; edge-connectivity; minimum degree

 \end{abstract}

\section{Introduction}

We use Bondy and Murty \cite{bon} for terminology and notations not
defined here and consider finite simple graphs only.

Let $G=(V(G),E(G))$ be a connected graph. We use $n(G), e(G),
\kappa(G), \kappa'(G)$ and $\delta(G)$ to denote the \emph{order,
size, connectivity, edge-connectivity} and \emph{minimum degree}
of $G$, respectively. Let $u$ be a vertex of $G$. We use $N_{G}(u)$
to denote the set of vertices which is adjacent with $u$ (also
called the \emph{neighbors} of $u$) in the graph $G$. Let $S$ be a
subset of $V(G)$(or $E(G)$). The \emph{induced subgraph} of $G$ is
denoted by $G[S]$. Furthermore, we use $G-S$ to denote the subgraph
$G[V(G)\backslash S]$(or $G[E(G)\backslash S]$), respectively.
For $x, y \in V(G)$, the length of a shortest path joining $x$ and $y$ is called the \emph{distance} between $x$ and
$y$ and denoted by $d_{G}(x,y)$. The \emph{diameter} of a graph $G$, denoted by $dim(G)$, is the greatest distance between two vertices of $G$.

Let $H$ be a given graph. A graph $G$ is said to be \emph{$H$-free}
if $G$ contains no induced copies of $H$. If $G$ is $H$-free, then
$H$ is called a \emph{forbidden subgraph} of $G$. Note that if $H_{1}$
is an induced subgraph of $H_{2}$, then every $H_{1}$-free graph is
also $H_{2}$-free. For a class of graphs $\mathcal{H}$, the graph $G$ is \emph{$\mathcal{H}$-free} if
$G$ is $H$-free for every $H \in \mathcal{H}$. For two sets $\mathcal{H}_{1}$ and $\mathcal{H}_{2}$ of connected graphs,
we write $\mathcal{H}_{1} \preceq \mathcal{H}_{2}$ if for every graph $H_{2} \in \mathcal{H}_{2}$,
there exists a graph $H_{1} \in \mathcal{H}_{1}$ such that $H_{1}$ is an induced subgraph of $H_{2}$.
If $\mathcal{H}_{1} \preceq \mathcal{H}_{2}$, then every $\mathcal{H}_{1}$-free graph is
also $\mathcal{H}_{2}$-free.

\begin{figure}[htp]
\begin{centering}
\includegraphics[width=12cm]{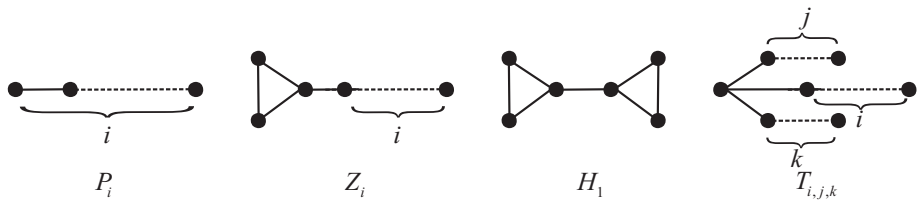}
\par\end{centering}
\caption{Some special graphs: $P_{i}, Z_{i}, H_{1}$ and $T_{i,j,k}$.}
\label{f1}
\end{figure}

As usual, we use $K_{n}$ to denote the complete graph of order $n$,
and $K_{m,n}$ to denote the complete bipartite graph with partition
sets of size $m$ and $n$. So the $K_{1}$ is a vertex, $K_{3}$ is a
triangle, $K_{1,r}$ is a star (the $K_{1,3}$ is also called a claw).
The clique $C$ is a subgraph of a graph $G$ such that
$G[V(C)]$ is a complete graph, and the \emph{clique number} $\omega(G)$ of a graph $G$ is the maximum
cardinality of a clique of $G$. Then we will show some special graphs
which are needed: (see Figure~\ref{f1})
\begin{itemize}
  \item $P_{i}$, the path with $i$ vertices (note that $P_{1} = K_{1}$ and $P_{2} = K_{2}$);
  \item $Z_{i}$, a graph obtained by identifying a vertex of a $K_{3}$ with an end-vertex of a $P_{i+1}$;
  \item $H_{1}$, a graph obtained by identifying a vertex of a $K_{3}$ with the one-degree vertex of a $Z_{1}$;
  \item $T_{i,j,k}$, a graph consisting of three paths $P_{i + 1}, P_{j + 1}$ and $P_{k + 1}$ with the common starting vertex.
\end{itemize}

Let $X$ and $Y$ be nonempty subsets of $V(G)$, we denote by $E[X, Y]$ the set of edges of $G$
with one end in $X$ and the other end in $Y$, and by $e(X,Y)$ their number. When $Y = V(G)\backslash X$,
the set $E[X, Y]$ is called the \emph{edge cut} of $G$ associated with $X$.
The edge cut set $S$ with the minimum number of edges is called the \emph{minimum edge cut}. It is well-known that $\kappa(G) \leq \kappa'(G) \leq \delta(G)$. In \cite{wang}, Wang, Tsuchiya and Xiong characterize all the pairs $R, S$ such that every connected $\{R, S\}$-free graph $G$ has $\kappa(G) = \kappa'(G)$.

\begin{theorem} \label{wang1} \emph{(Wang, Tsuchiya and Xiong \cite{wang})}
Let $S$ be a connected graph. Then $G$ being a connected $S$-free graph implies $\kappa(G) = \kappa'(G)$ if and only if $S$ is an induced subgraph of $P_{3}$.
\end{theorem}

\begin{theorem} \label{wang2} \emph{(Wang, Tsuchiya and Xiong \cite{wang})}
Let $\mathcal{H} = \{R, S\}$ be a set of two connected graphs such
that $R, S \neq P_{3}$. Then $G$ being a connected
$\mathcal{H}$-free graph implies $\kappa(G) = \kappa'(G)$ if and
only if $\mathcal{H} \preceq \{Z_{1}, P_{5}\}$, $\mathcal{H} \preceq
\{Z_{1}, K_{1, 4}\}$, $\mathcal{H} \preceq \{Z_{1}, T_{1, 1, 2}\}$,
$\mathcal{H} \preceq \{P_{4}, H_{0}\}$ or $\mathcal{H} \preceq
\{K_{1, 3}, H_{0}\}$, where $H_{0}$ is the unique simple graph with
degree sequence 4, 2, 2, 2, 2.
\end{theorem}

In \cite{hel}, Hellwig and Volkmann introduce a lot of sufficient conditions for $\kappa'(G) = \delta(G)$.

\begin{theorem} \label{hel}
Let $G$ be a connected graph satisfying the one of following
conditions:
\begin{enumerate}
  \item \emph{(Chartrand \cite{cha})} $n(G) \leq 2\delta(G) + 1$,
  \item \emph{(Lesniak \cite{les})} $d_{G}(u) + d_{G}(v) \geq n(G) - 1$ for all pairs $u, v$ of nonadjacent vertices,
  \item \emph{(Plesn\'{\i}k \cite{ple})} $dim(G) = 2$,
  \item \emph{(Volkmann \cite{vol})} $G$ is bipartite and $n(G) \leq 4\delta(G) - 1$,
  \item \emph{(Plesn\'{\i}k and Zn\'am \cite{ple1})} there are no four vertices $u_{1}, u_{2}, v_{1}, v_{2}$ with\\ $d_{G}(u_{1}, u_{2}), d_{G}(u_{1}, v_{2}), d_{G}(v_{1}, u_{2}), d_{G}(v_{1}, v_{2}) \geq 3$,
  \item \emph{(Plesn\'{\i}k and Zn\'am \cite{ple1})} $G$ is bipartite and $dim(G) = 3$,
  \item \emph{(Xu \cite{xu})} there exist $\lfloor n(G)/2 \rfloor$ pairs $(u_{i}, v_{i})$ of vertices such that $d_{G}(u_{i}) + d_{G}(v_{i}) \geq n(G)$ for $i = 1, 2, \cdots, \lfloor n(G)/2 \rfloor$,
  \item \emph{(Dankelmann and Volkmann \cite{dan})} $\omega(G) \leq p$ and $n(G) \leq 2\lfloor p\delta(G)/(p - 1) \rfloor - 1$.
\end{enumerate}
Then $\kappa'(G) = \delta(G)$.
\end{theorem}

In this paper, we also consider and characterize the forbidden subgraphs for
$\kappa'(G) = \delta(G)$.

\begin{theorem} \label{main1}
Let $S$ be a connected graph. Then $G$ being a connected $S$-free graph implies $\kappa'(G) = \delta(G)$ if and only if $S$ is an induced subgraph of $P_{4}$.
\end{theorem}

\begin{theorem}\label{main2}
Let $\mathcal{H} = \{R, S\}$ be a set of two connected graphs such
that $R$ and $S$ are not an induced subgraph of $P_{4}$. Then $G$
being a connected $\mathcal{H}$-free graph implies $\kappa'(G) =
\delta(G)$ if and only if $\mathcal{H} \preceq \{H_{1}, P_{5}\}$,
$\mathcal{H} \preceq \{Z_{2}, P_{6}\}$, or $\mathcal{H} \preceq
\{Z_{2}, T_{1, 1, 3}\}$.
\end{theorem}

Note that all families of connected graphs satisfies $\kappa(G) <
\kappa'(G)$ or $\kappa'(G) < \delta(G)$ should be $\kappa(G) <
\delta(G)$. By Theorems \ref{wang1} and \ref{wang2}, we may get the
following corollaries.

\begin{corollary} \label{co1}
Let $S$ be a connected graph. Then $G$ being a connected $S$-free
graph implies $\kappa(G) = \delta(G)$ if and only if $S$ is an
induced subgraph of $P_{3}$.
\end{corollary}

\begin{corollary} \label{co2}
Let $\mathcal{H} = \{R, S\}$ be a set of two connected graphs such
that $R$ and $S$ are not an induced subgraph of $P_{3}$. Then $G$
being a connected $\mathcal{H}$-free graph implies $\kappa(G) =
\delta(G)$ if and only if $\mathcal{H} \preceq \{H_{0}, P_{4}\}$,
$\mathcal{H} \preceq \{Z_{1}, P_{5}\}$, or $\mathcal{H} \preceq \{Z_{1}, T_{1, 1,
2}\}$.
\end{corollary}

In fact, we also present a general result as follow. Now
Corollaries \ref{co1} and \ref{co2} follow easily from  Theorems~\ref{wang1}, \ref{wang2}, \ref{main1}, \ref{main2} and \ref{c}. Note that $P_4$ may be one of  the pair of forbidden subgraphs, see Theorem~\ref{main2}, while $P_4$ is the forbidden subgraph from Theorem~\ref{main1}, this means that the other subgraph may be any subgraph of $G$ when $P_4$ is one of  a pair of forbidden subgraphs.

\begin{theorem}\label{c}
Let $G$ be a connected graph, and $f(G), g(G), t(G)$ are three
invariants of $G$ with $f(G) \leq g(G) \leq t(G)$. If the following
statements hold:
\begin{enumerate}
  \item $G$ is $\mathcal{H}$-free implies $f(G) = g(G)$ if and only
  if $\mathcal{H} \in \mathbf{H_{1}}$;
  \item $G$ is $\mathcal{H}$-free implies $g(G) = t(G)$ if and only
  if $\mathcal{H} \in \mathbf{H_{2}}$,
\end{enumerate}
then $G$ is $\mathcal{H}$-free implies $f(G) = t(G)$ if and only
  if $\mathcal{H} \in \mathbf{H_{1}} \bigcap \mathbf{H_{2}}$.
  Here $\mathbf{H_{i}}$ is the set of class of given graphs, i.e.,
each element of $\mathbf{H_{i}}$ is a class of given graphs
$\mathcal{H}$, for $i \in \{1, 2\}$. $\mathbf{H_{1}} \bigcap
\mathbf{H_{2}} := \{\mathcal{H}_{1} \bigcap \mathcal{H}_{2}|
\mathcal{H}_{1} \in \mathbf{H_{1}}, \mathcal{H}_{2} \in
\mathbf{H_{2}}$, and $|\mathcal{H}_{1}| = |\mathcal{H}_{2}|\}$, and
$\mathcal{H}_{1} \bigcap \mathcal{H}_{2}$ is the set with order
$|\mathcal{H}_{1}|$, which each element is the common induced
subgraph of one graph in $\mathcal{H}_{1}$ and one graph in
$\mathcal{H}_{2}$, respectively.
\end{theorem}

{\bf Proof.}
First suppose $G$ is $\mathcal{H}$-free and $\mathcal{H} \in
\mathbf{H_{1}} \bigcap \mathbf{H_{2}}$, then $\mathcal{H} \in
\mathbf{H_{1}}$ and $\mathcal{H} \in \mathbf{H_{2}}$. By (1) and
(2), $f(G) = g(G)$ and $g(G) = t(G)$. It means that $f(G) = g(G) =
t(G)$. This completes the sufficiency.

Now we prove the necessity. Suppose that $f(G)=t(G)$. Then $f(G)=g(G)=t(G)$ since $f(G) \leq g(G) \leq t(G)$. Therefore, both $\mathcal{H} \in \mathbf{H_{1}}$ and $\mathcal{H}
\in \mathbf{H_{2}}$ must  hold, by (1) and (2).   It means that $\mathcal{H} \in
\mathbf{H_{1}} \bigcap \mathbf{H_{2}}$. This completes the
proof.\qed

\section{The necessity part of Theorems \ref{main1} and \ref{main2}}

\begin{figure}[H]
\begin{centering}
\includegraphics[width=14cm]{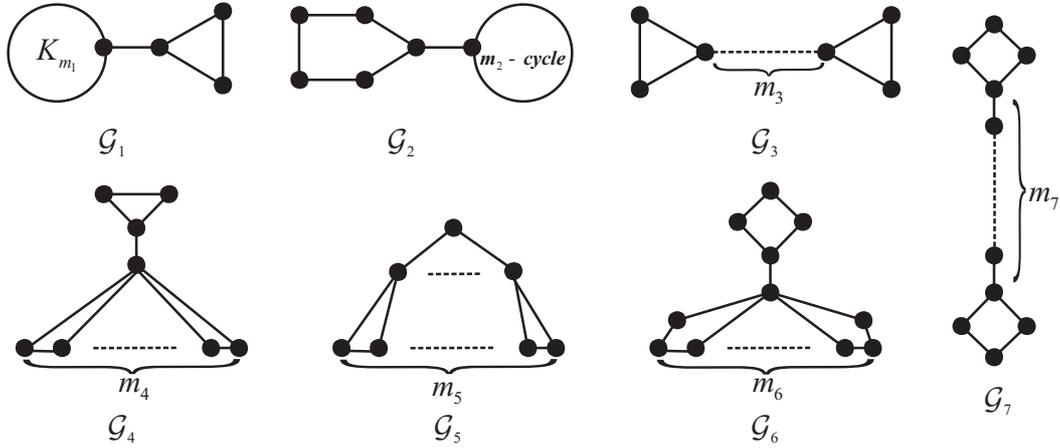}
\par\end{centering}
\caption{Some classes of graphs satisfies $1 = \kappa'(G) <
\delta(G) = 2$.}
\label{f2}
\end{figure}

We first construct some families of connected graphs
$\mathcal{G}_{i}, i=1,\cdots,7$ (see Figure~\ref{f2}). It is easy to
see that each $G \in \mathcal{G}_{i}$ satisfies $1 = \kappa'(G) <
\delta(G) = 2$.

\textbf{The necessity part of Theorem \ref{main1}.} Let $S$ be a
graph such that every connected $S$-free graph is $\kappa'(G) =
\delta(G)$. Then $S$ is an induced subgraph of all graphs in
$\mathcal{G}_{i}, i=1,\cdots,7$.

Note that the common induced
subgraph of the graphs in $\mathcal{G}_{1}$ and $\mathcal{G}_{2}$ is
a path. Since the largest induced path of the graphs in
$\mathcal{G}_{1}$ is $P_{4}$, $S$ must be an induced subgraph of
$P_{4}$. This completes the proof of the necessity part of Theorem
\ref{main1}.\qed

\textbf{The sufficiency part of Theorem \ref{main2}.} Let $R$ and
$S$ are not an induced subgraph of $P_{4}$ graphs such that every
connected $\{R, S\}$-free graph is $\kappa'(G) = \delta(G)$. Then
all graphs in $\mathcal{G}_{i}, i=1,\cdots,7$ should contain either
$R$ or $S$ as an induced subgraph. Without loss of generality, we
may assume that $R$ is an induced subgraph of all graphs in
$\mathcal{G}_{1}$. Note that all graphs in $\mathcal{G}_{1}$ contain no induced cycle with length at least 4 as an induced
subgraph, so we need to consider the following four cases.

\textbf{Case 1.} $R$ contain a clique $K_{t}$ with $t \geq 4$.

It means that, for $i \in \{2, 3, 4, 5, 6, 7\}$, all graphs in
$\mathcal{G}_{i}$ are $R$-free, and should contain $S$ as an induced
subgraph. Note that all graphs in $\mathcal{G}_{2}$ are
$K_{3}$-free, and all graphs in $\mathcal{G}_{3}$ are $K_{1,
3}$-free, so $S$ should be a path. Since the largest induced path of
the graphs in $\mathcal{G}_{4}$ is $P_{4}$, $S$ should be an induced
subgraph of $P_{4}$, a contradiction.

\textbf{Case 2.} $R$ don't contain the clique $K_{t}$ with $t \geq
4$, but contain two $K_{3}$.

Since $R$ is an induced subgraph of all graphs in $\mathcal{G}_{1}$,
$R$ should be $H_{1}$. It means that, for $i \in \{2, 3, 5, 6, 7\}$,
all graphs in $\mathcal{G}_{i}$ are $R$-free, and should contain $S$
as an induced subgraph. Note that all graphs in $\mathcal{G}_{2}$
are $K_{3}$-free, and all graphs in $\mathcal{G}_{3}$ are $K_{1,
3}$-free, so $S$ should be a path. Since the largest induced path of
the graphs in $\mathcal{G}_{5}$ is $P_{5}$, $S$ should be an induced
subgraph of $P_{5}$. So $\mathcal{H} = \{R, S\} \preceq \{H_{1},
P_{5}\}$.

\textbf{Case 3.} $R$ don't contain the clique $K_{t}$ with $t \geq
4$, but contain exactly one $K_{3}$.

Since $R$ is an induced subgraph of all graphs in $\mathcal{G}_{1}$,
$R$ should be an induced subgraph of $Z_{2}$. It means that, for $i
\in \{2, 6, 7\}$, all graphs in $\mathcal{G}_{i}$ are $R$-free, and
should contain $S$ as an induced subgraph. Note that the common
induced subgraph of all graphs in $\mathcal{G}_{2}$ and
$\mathcal{G}_{7}$ are a tree with the maximum degree 3 or a path. If
$S$ is a tree with the maximum degree 3, since the common induced
tree with the maximum degree 3 of all graphs in $\mathcal{G}_{6}$
and $\mathcal{G}_{7}$ are $T_{1, 1, 3}$, $S$ should be an induced
subgraph of $T_{1, 1, 3}$. So $\mathcal{H} = \{R, S\} \preceq
\{Z_{2}, T_{1, 1, 3}\}$. If $S$ is a path. Since the largest induced
path of the graphs in $\mathcal{G}_{6}$ is $P_{6}$, $S$ should be an
induced subgraph of $P_{6}$. So $\mathcal{H} = \{R, S\} \preceq
\{Z_{2}, P_{6}\}$.

\textbf{Case 4.} $R$ is a tree.

Since all graphs in $\mathcal{G}_{1}$ are $K_{1, 3}$-free, $R$
should be a path. Note that the largest induced path of the graphs
in $\mathcal{G}_{1}$ is $P_{4}$, so $R$ should be an induced
subgraph of $P_{4}$, a contradiction.

From the proofs above, we have that $\mathcal{H} \preceq \{H_{1},
P_{5}\}$, $\mathcal{H} \preceq \{Z_{2}, P_{6}\}$, or $\mathcal{H}
\preceq \{Z_{2}, T_{1, 1, 3}\}$. This completes the proof of the
necessity part of Theorem \ref{main2}.\qed

\section{The sufficiency part of Theorems \ref{main1} and \ref{main2}}
\hspace*{\parindent}

\textbf{The sufficiency part of Theorem \ref{main1}.} Let $G$ be a
connected $P_{4}$-free graph. Then $dim(G) \leq 2$. If $dim(G) = 1$,
$G$ must be a complete graph and $\kappa'(G) = \delta(G) = n - 1$.
If $dim(G) = 2$, by Theorem \ref{hel} (3), $\kappa'(G) = \delta(G)$.
This completes the proof of the sufficiency part of Theorem
\ref{main1}.\qed

\textbf{The sufficiency part of Theorem \ref{main2}.} Let $G$ be a
connected $\mathcal{H}$-free graph such that $\kappa'(G) <
\delta(G)$, where $\mathcal{H} \preceq \{H_{1}, P_{5}\}, \{Z_{2},
P_{6}\}$, or $\{Z_{2}, T_{1, 1, 3}\}$. Then there must exists a
minimum edge cut, say $M$, such that $|M| = \kappa'(G) < \delta(G)$.
Let $G_{1}$ and $G_{2}$ are the components of $G - M$, and let
$S_{i} = V(G_{i}) \bigcap V(M)$, $i \in \{1, 2\}$. Then $|S_{i}|
\leq |M| = \kappa'(G) < \delta(G)$, say $|S_{i}| = s_{i}$, $i \in
\{1, 2\}$.

\begin{claim} \label{cla1}
For $i \in \{1, 2\}$, $V(G_{i} - S_{i}) \neq \emptyset$. Moreover,
for any $x \in V(G_{i} - S_{i})$, $N_{G}(x) \bigcap V(G_{i} - S_{i})
\neq \emptyset$.
\end{claim}

\textbf{Proof}. We will count the number of edges of $G_{i}$ for $i
\in \{1, 2\}$.
\begin{align*}
|E(G_{i})| & = \frac{1}{2}\left(\sum\limits_{x \in V(G_{i})}d_{G}(x)
- \kappa'(G)\right) \\
& \geq \frac{1}{2}\left(\delta(G) |V(G_{i})| - \kappa'(G)\right) \\
& \geq \frac{1}{2}\left(\delta(G) s_{i} -
\kappa'(G)\right) \\
& > \frac{1}{2}\kappa'(G)\left(s_{i} - 1 \right) \\
& \geq \frac{1}{2}s_{i}\left(s_{i} - 1 \right)
\end{align*}

Note that the complete graph $K_{s_{i}}$ has
$\frac{1}{2}s_{i}\left(s_{i} - 1 \right)$ edges. It means that
$|V(G_{i})| > s_{i}$, i.e., $V(G_{i} - S_{i}) \neq \emptyset$.

Moreover, for any $x \in V(G_{i} - S_{i})$, since $d_{G}(x) \geq
\delta(G) > \kappa'(G) \geq s_{i}$, $N_{G}(x) \bigcap V(G_{i} -
S_{i}) \neq \emptyset$. This completes the proof of Claim
\ref{cla1}. \qed

Now we will distinguish the following two cases to complete our
proof.

\textbf{Case 1.} $G$ contains a $P_{4} = x_{0}x_{1}x_{2}x_{3}$ with
$x_{0} \in V(G_{1} - S_{1}), x_{1} \in S_{1}, x_{2} \in S_{2}$, and
$x_{3} \in V(G_{2} - S_{2})$.

\textbf{Subcase 1.1.} $\mathcal{H} \preceq \{H_{1}, P_{5}\}$.

By Claim \ref{cla1}, there exist two vertices $x'_{0} \in V(G_{1} -
S_{1})$ and $x'_{3} \in V(G_{2} - S_{2})$ such that $x_{0}x'_{0},
x_{3}x'_{3} \in E(G)$. Then $G[\{x'_{0}, x_{0}, x_{1}, x_{2}, x_{3},
x'_{3}\}] \cong H_{1}$ (if  $x_{1}x'_{0}, x_{2}x'_{3} \in E(G)$),
or $G[\{x'_{0}, x_{0}, x_{1}, x_{2}, x_{3}\}] \cong P_{5}$ (if
$x_{1}x'_{0} \notin E(G)$), or $G[\{x_{0}, x_{1}, x_{2}, x_{3},
x'_{3}\}] \cong P_{5}$ (if  $x_{2}x'_{3} \notin E(G)$), a
contradiction.

\textbf{Subcase 1.2.} $\mathcal{H} \preceq \{Z_{2}, P_{6}\}$.

By Claim \ref{cla1}, there exist two vertices $x'_{0} \in V(G_{1} -
S_{1})$ and $x'_{3} \in V(G_{2} - S_{2})$ such that $x_{0}x'_{0},
x_{3}x'_{3} \in E(G)$. Then $G[\{x'_{0}, x_{0}, x_{1}, x_{2}, x_{3},
x'_{3}\}] \cong P_{6}$ (if  $x_{1}x'_{0}, x_{2}x'_{3} \notin
E(G)$), or $G[\{x'_{0}, x_{0}, x_{1}, x_{2}, x_{3}\}] \cong Z_{2}$
(if  $x_{1}x'_{0} \in E(G)$), or $G[\{x_{0}, x_{1}, x_{2}, x_{3},
x'_{3}\}] \cong Z_{2}$ (if  $x_{2}x'_{3} \in E(G)$), a
contradiction.

\textbf{Subcase 1.3.} $\mathcal{H} \preceq \{Z_{2}, T_{1, 1, 3}\}$.

By Claim \ref{cla1}, $N_{G}(x_{0}) \bigcap V(G_{1} - S_{1}) \neq
\emptyset$ and $N_{G}(x_{3}) \bigcap V(G_{2} - S_{2}) \neq
\emptyset$.

Suppose that  $|N_{G}(x_{0}) \bigcap V(G_{1} - S_{1})| \geq 2$ or
$|N_{G}(x_{3}) \bigcap V(G_{2} - S_{2})| \geq 2$. Without loss of
generality, we may assume that $|N_{G}(x_{0}) \bigcap V(G_{1} -
S_{1})| \geq 2$, it means there exist two vertices $x'_{0}, x''_{0}
\in V(G_{1} - S_{1})$ such that $x_{0}x'_{0}, x_{0}x''_{0} \in
E(G)$. Then $G[\{x'_{0}, x''_{0}, x_{0}, x_{1}, x_{2}, x_{3}\}]
\cong T_{1, 1, 3}$ (if  $x'_{0}x''_{0}, x'_{0}x_{1}, x''_{0}x_{1}
\notin E(G)$), or $G[\{x'_{0}, x''_{0}, x_{0}, x_{1}, x_{2}\}] \cong
Z_{2}$ (if  $x''_{0}x'_{0} \in E(G)$ and $x'_{0}x_{1},
x''_{0}x_{1}\notin E(G)$), or $G[\{x'_{0}, x_{0}, x_{1}, x_{2},
x_{3}\}] \cong Z_{2}$ (if  $x'_{0}x_{1} \in E(G)$), or
$G[\{x''_{0}, x_{0}, x_{1}, x_{2}, \\  x_{3}\}] \cong Z_{2}$ (if
$x''_{0}x_{1} \in E(G)$), a contradiction.

Suppose that  $N_{G}(x_{0}) \bigcap V(G_{1} - S_{1}) = \{x'_{0}\}$ and
$N_{G}(x_{3}) \bigcap V(G_{2} - S_{2}) = \{x'_{3}\}$. Note that
$N_{G}(x_{0}) \subseteq \{x'_{0}\} \bigcup S_{1}$ and $N_{G}(x_{3})
\subseteq \{x'_{3}\} \bigcup S_{2}$. Then $d_{G}(x_{0}) \leq s_{1} +
1$ and $d_{G}(x_{3}) \leq s_{2} + 1$. Since $d_{G}(x_{0}) \geq
\delta(G) > \kappa'(G) \geq s_{1}$ and $d_{G}(x_{3}) \geq \delta(G)
> \kappa'(G) \geq s_{2}$, $d_{G}(x_{0}) \geq s_{1} + 1$ and $d_{G}(x_{3}) \geq s_{2} + 1$.
Therefore  $d_{G}(x_{0}) = s_{1} + 1$ and $d_{G}(x_{3}) = s_{2} +
1$. It means that $N_{G}(x_{0}) = S_{1} \bigcup \{x'_{0}\}$,
$N_{G}(x_{3}) = S_{2} \bigcup \{x'_{3}\}$, and $s_{1} = s_{2} =
\kappa'(G)$. Since $|M| = \kappa'(G) = s_{1} = s_{2}$, the each
vertex in $S_{1}$ is just adjacent to exactly one vertex which is in
$S_{2}$, and vice versa. Suppose $s_{1} \geq 2$. Then there exists a
vertex $x'_{1} \in S_{1}$ such that $x'_{1} \neq x_{1}$. Therefore
$G[\{x'_{0}, x_{0}, x'_{1}, x_{1}, x_{2}, x_{3}\}] \cong T_{1, 1,
3}$ (if  $x'_{0}x'_{1}, x'_{0}x_{1}, x'_{1}x_{1} \notin E(G)$), or
$G[\{x'_{0}, x'_{1}, x_{0}, x_{1}, x_{2}\}] \cong Z_{2}$ (if
$x'_{0}x'_{1} \in E(G)$ and $x'_{0}x_{1}, x'_{1}x_{1}\notin E(G)$),
or $G[\{x'_{0}, x_{0}, x_{1}, x_{2}, x_{3}\}] \cong Z_{2}$ (if
$x'_{0}x_{1} \in E(G)$), or $G[\{x'_{1}, x_{0}, x_{1}, x_{2},
x_{3}\}] \cong Z_{2}$ (if  $x'_{1}x_{1} \in E(G)$), a
contradiction. Suppose $s_{1} = 1$. Then $s_{2} = \kappa'(G) = 1$
and $\delta(G) = 2$. Assume $d_{G}(x_{1}) \geq 3$. Then there exists
a vertex $x'_{1} \in V(G_{1} - S_{1})$, such that $x'_{1}x_{1} \in
E(G)$ and $x'_{1} \neq x_{0}$. Therefore $G[\{x_{0}, x'_{1}, x_{1},
x_{2}, x_{3}, x'_{3}\}] \cong T_{1, 1, 3}$ (if  $x_{0}x'_{1},
x'_{3}x_{2} \notin E(G)$), or $G[\{x_{0}, x'_{1}, x_{1}, x_{2},
x_{3}\}] \cong Z_{2}$ (if  $x_{0}x'_{1} \in E(G)$), or $G[\{x'_{3},
x_{3}, x_{2}, x_{1}, x_{0}\}] \\ \cong Z_{2}$ (if  $x'_{3}x_{2} \in
E(G)$), a contradiction. Assume $d_{G}(x_{1}) = 2$. Then it means
that $N_{G}(x_{1}) = \{x_{0}, x_{2}\}$ and $d_{G}(x) = d_{G_{1}}(x)$
for any $x \in V(G_{1} - \{x_{0}, x_{1}\})$. Since $\delta(G) = 2$
and $d_{G_{1} - S_{1}}(x_{0}) = 1$, there exist some vertices in
$V(G_{1} - S_{1})$ such that their degree in $G$ are at least 3.
Then we choose a vertex $y_{0} \in V(G_{1} - S_{1})$, such that
$d_{G}(y_{0}) \geq 3$ and $d_{G}(y_{0}, x_{1})$ as small as
possible. Let $P'$ is the shortest path between $x_{1}$ and $y_{0}$.
Then all inner vertices of $P'$ should have degree two. Let $y_{1},
y_{2} \in N_{G}(y)$ and $y_{1}, y_{2} \notin V(P')$. Then
$G[\{y_{1}, y_{2}, x_{2}, x_{3}\} \bigcup V(P')]$ contians an
induced $T_{1, 1, 3}$ (if  $y_{1}y_{2} \notin E(G)$), or
$G[\{y_{1}, y_{2}, x_{2}\} \bigcup V(P')]$ contians an induced $
Z_{2}$ (if  $y_{1}y_{2} \in E(G)$), a contradiction.

\textbf{Case 2.} $G$ contains no $P_{4} = x_{0}x_{1}x_{2}x_{3}$ with
$x_{0} \in V(G_{1} - S_{1}), x_{1} \in S_{1}, x_{2} \in S_{2}$, and
$x_{3} \in V(G_{2} - S_{2})$.

Let $S^{1}_{i}=\{x \in S_{i}: N_{G}(x) \bigcap V(G_{i} - S_{i}) \neq
\emptyset\}$, and $S^{2}_{i} = S_{i} - S^{1}_{i}$ for $i = 1, 2$.
Then $S^{2}_{i} \neq \emptyset$ and $E(S^{1}_{1}, S^{1}_{2}) =
\emptyset$. By the minimality of $M$ and the definition of $S_{i}$,
$E(S^{1}_{i}, S^{2}_{i}), E(S^{1}_{1}, S^{2}_{2}), E(S^{2}_{1},
S^{1}_{2}) \neq \emptyset$. Now we choose a path $P_{0}$ between
$x_{1}$ and $x_{2}$, such that $x_{1} \in S^{1}_{1}$ and $x_{2} \in
S^{1}_{2}$, and the length of path as small as possible. Then
$|V(P_{0})| \geq 3$ and all inner vertices of $P_{0}$ must be in
$S^{2}_{i}$. Let $x_{0} \in V(G_{1} - S_{1})$ and $x_{3} \in V(G_{2}
- S_{2})$, such that $x_{0}x_{1}, x_{2}x_{3} \in E(G)$. Then
$G[V(P_{0}) \bigcup \{x_{0}, x_{3}\}]$ is an induced path with at
least 5 vertices, say $P_{1}$.

\textbf{Subcase 2.1.} $\mathcal{H} \preceq \{H_{1}, P_{5}\}$.

$P_{1}$ is an induced path with at least 5 vertices, a
contradiction.

\textbf{Subcase 2.2.} $\mathcal{H} \preceq \{Z_{2}, P_{6}\}$.

By Claim \ref{cla1}, there exist a vertex $x'_{0} \in V(G_{1} -
S_{1})$ such that $x_{0}x'_{0} \in E(G)$. Then $G[\{x'_{0}\} \bigcup
V(P_{1})]$ contians an induced $P_{6}$ (if  $x_{1}x'_{0} \notin
E(G)$), or an induced $Z_{2}$ (if  $x_{1}x'_{0} \in E(G)$), a
contradiction.

\textbf{Subcase 2.3.} $\mathcal{H} \preceq \{Z_{2}, T_{1, 1, 3}\}$.

By Claim \ref{cla1} and $|S^{1}_{1}| < s_{1} < \delta(G)$, there
exist two vertices $x'_{0}, x''_{0} \in V(G_{1} - S_{1})$ such that
$x_{0}x'_{0}, x_{0}x''_{0} \in E(G)$. Then $G[\{x'_{0}, x''_{0}\}
\bigcup V(P_{1})]$ contians an induced $T_{1, 1, 3}$ (if
$x'_{0}x''_{0}, x'_{0}x_{1}, x''_{0}x_{1} \notin E(G)$), or an
induced $Z_{2}$ (if  $x''_{0}x'_{0} \in E(G)$ and $x'_{0}x_{1},
x''_{0}x_{1}\notin E(G)$), or an induced $Z_{2}$ (if  $x'_{0}x_{1}
\in E(G)$ or $x''_{0}x_{1} \in E(G)$), a contradiction.

This completes the proof of the sufficiency part of Theorem
\ref{main2}.\qed

\section{Concluding remark}

In this paper, we give a completely characterzation of all pairs
$\{R,S\}$ of graphs such that every connected $\{R,S\}$-free graph
has the same edge-connectivity and minimum degree. All graphs in
Figure 2 have edge-connectivity one, we also can construct some
graphs for arbitrarily large edge-connectivity to show that Theorem
\ref{main1} also hold. But for forbidden pairs $\mathcal{H} = \{R,
S\}$, we have not enough graphs to see that whether could get more
wide forbidden pairs to guarantee the graphs having the same
edge-connectivity and minimum degree, when we increase the
edge-connectivity.


\begin{thebibliography}{20}

\bibitem{bon} J. A. Bondy and U. S. R. Murty, Graph theory with applications, Macmillan, London and Elsevier, New York, 1976.

\bibitem{cha} G. Chartrand, A graph-theoretic approach to a communications problem, SIAM J. Appl. Math., 14 (1966) 778-781.

\bibitem{dan} P. Dankelmann and L. Volkmann, New sufficient conditions for equality of minimum degree and edge-connectivity, Ars Combin., 40 (1995) 270-278.

\bibitem{hel} A. Hellwig and L. Volkmann, Sufficient conditions for graphs to be $\lambda'$-optimal, super-edge-connected, and maximally edge-connected, J. Graph Theory, 48 (2005) 228-246.

\bibitem{les} L. Lesniak, Results on the edge-connectivity of graphs, Discrete Math., 8 (1974) 351-354.

\bibitem{ple} J. Plesn\'{\i}k, Critical graphs of given diameter, Acta Fac. Rerum Nat. Univ. Commenian. Math., 30 (1975) 71-93.

\bibitem{ple1} J. Plesn\'{\i}k and S. Zn\'am, On equality of edge-connectivity and minimum degree of a graph, Arch. Match.(Brno), 25 (1989)  19-25.

\bibitem{vol} L. Volkmann, Bemerkungen zum $p$-fachen Kantenzusammenhang von Graphen, An. Univ. Bucuresti Mat., 37 (1988) 75-79.

\bibitem{wang} S. Wang, S. Tsuchiya and L. Xiong, Forbidden pairs for equality of connectivity and edge-connectivity of graphs, submitted.

\bibitem{xu} J. -M. Xu, A sufficient condition for equality of arc-connectivity and minimum degree of a digraph, Discrete Math., 133 (1994) 315-318.

\end{thebibliography}
\end{document}